\documentclass[11pt]{amsart}
\usepackage{amsmath,amssymb,mathrsfs}
\usepackage{graphicx}

\title{Equivariant deformations of LeBrun's self-dual metrics with torus action}
\author{Nobuhiro Honda 
}  

\newcommand{\ol}{\overline}
\newcommand{\ra}{\rightarrow}
\newcommand{\lra}{\longrightarrow}

\newcommand{\set}{\,|\,}
\newcommand{\proofend}{\hfill$\square$}

\newcommand{\vsp}{\vspace{3mm}}
\newcommand{\ptl}{\partial}

\newtheorem{prop}{Proposition}[section]

\newtheorem{example}[prop]{Example}

\begin{document}

\maketitle

\begin{abstract} 
We investigate $U(1)$-equivariant deformations of C. LeBrun's self-dual metric with torus action. 
We explicitly determine all $U(1)$-subgroups of the torus for which one can obtain  $U(1)$-equivariant deformation that do not preserve semi-free $U(1)$-action.
This gives many new self-dual metrics with $U(1)$-action which are not conformally isometric to LeBrun metric.
We also count the dimension of the moduli space of  self-dual metrics with $U(1)$-action obtained in this way.
\end{abstract}



\section{Introduction} 
In \cite{LB91} C. LeBrun explicitly constructed a family of self-dual metrics on $n\mathbf{CP}^2$, the connected sum of $n$ copies of complex projective planes, where $n$ is an arbitrary positive integer.
His construction starts from giving distinct $n$ points on the upper half-space $\mathscr H^3$ with the usual hyperbolic metric.
Once these $n$ points are given, everything proceed in a canonical way.
Namely a principal  $U(1)$-bundle over the punctured $\mathscr H^3$ together with a connection is canonically constructed, and then on the total space of this $U(1)$-bundle a self-dual metric is naturally and explicitly introduced, for which the $U(1)$-action becomes isometric.
Then finally by choosing an appropriate conformal gauge (which is also concretely given), the self-dual metric is shown to extend to a compactification, yielding desired self-dual metric on $n\mathbf{CP}^2$.
Thus  LeBrun metrics on $n\mathbf{CP}^2$ are naturally parametrized by the set of different $n$ points on $\mathscr H^3$.

If the $n$ points are located in a general position, the corresponding LeBrun metric admits only a $U(1)$-isometry (coming from the principal bundle structure).
However, when the $n$ points are put in a collinear position, meaning that the $n$ points lie on the same geodesic on the hyperbolic $\mathscr H^3$, then the rotations around the geodesic can be lifted to the total space and it gives another $U(1)$-isometries of the LeBrun metric.
We call this kind of self-dual metrics on $n\mathbf{CP}^2$ LeBrun metric with torus action.
By a characterization theorem of LeBrun \cite{LB93}, being LeBrun metric with torus action is preserved under deformation keeping the torus action.

In this note, following a suggestion of LeBrun \cite[p.\,123, Remark]{LB93}, we investigate $U(1)$-isometric deformation of LeBrun metrics with torus action, where $U(1)$ is a subgroup of the torus.
In particular,  {\em{we determine all $U(1)$-subgroups of the torus for which one can obtain $U(1)$-equivariant deformation such that not full torus symmetry  survive.}}
Note that  on $2\mathbf{CP}^2$ every self-dual metric of positive scalar curvature is LeBrun metric with torus action \cite{P86} and such a subgroup in problem cannot exist  for $n=2$ (and also for $n=1$).
Of course, LeBrun's original $U(1)$-subgroup (coming from the principal bundle structure), which acts semi-freely on $n\mathbf{CP}^2$, has the desired property for $n\geq 3$.
We show that involving this subgroup, {\em{there are precisely $n$ numbers of $U(1)$-subgroups for which there exists the required equivariant deformation.}}
We concretely give these subgroups and observe that the remaining $(n-1)$ subgroups give non-LeBrun self-dual metric.
Also, we count the dimension of the moduli space of the resulting family of self-dual metrics with a non-semi-free $U(1)$-isometries.
Finally, we discuss some examples.


\section{Computation of the torus action on a cohomology group}
\noindent
\textbf{(2.1)}
Our proof of the main result is via twistor space.
So let $Z$ be the twistor space of a LeBrun metric  with torus action on $n\mathbf{CP}^2$.
In order to investigate $U(1)$-equivariant deformations of this metric, we calculate the torus action on the cohomology group $H^1(\Theta_Z)$ which is relevant to deformation of complex structure of $Z$.
In this subsection, to this end, we recall the explicit construction of $Z$  due to LeBrun \cite{LB91}.
We need to be careful in resolving singularities of a projective model of the twistor space, since in \cite{LB91} it is assumed that the semi-free $U(1)$-action does not extend to torus action, and since, under the existence of torus action, there are $n!$ possible ways of (small) resolutions and most of them do not yield a twistor space

First let $Q=\mathbf{CP}^1\times\mathbf{CP}^1$ be a quadratic surface and $\mathscr E\ra Q$ a rank-3 vector bundle
$$\mathscr E=\mathscr O(n-1,1)\oplus \mathscr O(1,n-1)\oplus \mathscr O\ra Q,$$
where $\mathscr O(k,l)$ denotes the line bundle over $Q$ whose bidegree is $(k,l)$.
Let $(\xi_0,\xi_1)$ (resp. $(\eta_0,\eta_1)$) be a homogeneous coordinate on the first (resp. the second) factor of $Q$, and set $U_0=\{(\xi_0,\xi_1)\set\xi_0\neq 0\}, V_0=\{(\eta_0,\eta_1)\set \eta_0\neq 0\}$.
On $U_0$ (resp. $V_0$) we use a non-homogeneous coordinate $u=\xi_1/\xi_0$ (resp. $v=\eta_1/\eta_0$).
We choose a trivialization of $\mathscr E$ over $U_0\times V_0\subset Q$, and let $(x,y,z)$ be the resulting fiber coordinate on $\mathscr E|_{U_0\times V_0}$.
Thus on the total space of $\mathscr E|_{U_0\times V_0}$ we can use $(u,v,x,y,z)$ as a global coordinate.
Then let $X$ be a compact (or complete) algebraic variety in $\mathbf P(\mathscr E)$ define by 
\begin{equation}\label{eqn-lb}
xy=z^2\prod_{i=1}^n (v-a_iu),
\end{equation}
where $a_1,a_2,\cdots,a_n$ are positive real numbers satisfying $a_1<a_2\cdots<a_n$.
((\ref{eqn-lb}) is an equation on $\mathbf P(\mathscr E|_{U_0\times V_0})$, but it can be naturally compactified in $\mathbf P(\mathscr E)$).
 $X$ has an obvious conic bundle structure over $Q$ whose discriminant locus is $C_1\cup C_2\cup\cdots C_n$, where $C_i$ is a $(1,1)$-curve in $Q$ defined by $v=a_iu$.
Further, the point $(x,y,z)=(0,0,1) \in\mathbf P(\mathscr E)$ lying over the fiber over the point $(u,v)=(0,0)$ is so called a compound $A_{n-1}$-singularity of $X$.
Similarly, by the choice of the degree of the direct summand in $\mathscr E$, the point $(0,0,1)\in\mathbf P(\mathscr E)$ over $(u,v)=(\infty,\infty)$ is also a compound  $A_{n-1}$-singularity of $X$.
We denote these two singularities of $X$ by $p_0$ and $\ol{p}_0$.
These are all the singularities of $X$.

We have to define a real structure.
In terms of the above coordinate $(u,v,x,y,z)$ on $\mathbf P(\mathscr E|_{U_0\times V_0})$ it is defined by 
\begin{equation}\label{eqn-rs}
\sigma:(u,v\,;x,y,z)\mapsto \left(\frac{1}{\ol v},\frac{1}{\ol u}\,;\frac{\ol y}{\ol u^{n-1}\ol v},\,\frac{\ol x}{\ol u\,\ol v^{n-1}},\, \ol z\right),
\end{equation} 
which preserves $X$, and interchanges the two singular points $p_0$ and $\ol{p}_0$ of $X$.

Next we give a small resolution of $p_0$.
To give it explicitly we write $\tilde x=x/z$ and $\tilde y=y/z$.
Then in an affine neighborhood of $p_0$ in $\mathbf P(\mathscr E)$, $X$ is defined by $\tilde x\tilde y=\prod_{i=1}^n(v-a_iu)$.
The small resolution of $p_0$ is a composition of $(n-1)$ blowing-ups, where the center is 2-dimensional in each steps;
As the first step we take a blow-up of $X$ along $\tilde x=v-a_{1}u=0$, yielding a new space $X_1$ and a morphism $X_1\ra X$.
Since this center is contained in $X$,  the exceptional locus $E_1$ arises only over $p_0$ and it is isomorphic to $\mathbf{CP}^1$.
Introducing a new coordinate $\tilde x_1$ by $\tilde x=\tilde x_1(v-a_{_1}u)$ on $E_1$, the new space $X_1$ is locally defined by $\tilde x_1\tilde y=\prod_{i\ge 2}(v-a_iu)$, so that having a compound $A_{n-2}$-singularity at the origin.
The second step is to blow-up $X_1$ along $\tilde x_1=v-a_2u=0$, giving a new space $X_2$ with a compound $A_{n-3}$-singularity at the new origin.
After repeating   this process $(n-1)$ times, the singularity $p_0$ is resolved, and the exceptional locus is a string of $(n-1)$ smooth rational curves.
This is how to obtain a small resolution of $p_0$.
Once a resolution of $p_0$ is given, another singularity $\ol{p}_0$ is naturally resolved by reality.
Let $Y\ra X$ be the small resolution of $p_0$ and $\ol{p}_0$ obtained in this way.
($Y$ is non-singular.)

Obviously, other small resolutions of $p_0$ can be obtained for each permutation of $n$ letters $\{1,2,\cdots,n\}$.
But keeping in mind that we have assumed $a_1<a_2<\cdots<a_n$ and that the curve $x=y=v-a_iu$ ($1\leq i\leq n$), which is over a discriminant locus $C_i\subset Q$, has to be a twistor line over the isolated fixed point of the torus action on $n\mathbf{CP}^2$, it is easily seen that if we take the resolution associated to a permutation other than $\{1,2,\cdots,n-1,n\}$ (giving the small resolution above) and $\{n,n-1,\cdots, 2, 1\}$, then the resulting space does not become a twistor space even after the blowing-down process which will explained next.

Next we explain the final step for obtaining the twistor space.
The conic bundle $X\ra Q$ has two distinct sections $E=\{x=z=0\}$ and $\ol{E}=\{y=z=0\}$, which are conjugate of each other.
These sections are disjoint from $p_0$ and $\ol{p}_0$ and their normal bundles  in $X$ are $\mathscr O(-1,1-n)$ and $\mathscr O(1-n,-1)$ respectively.
Clearly the small resolution $Y\ra X$ does not have any effect around $E$ and $\ol{E}$, so that it does not change the normal bundles.
Hence (if $n>2$) both $E$ and $\ol{E}$ (considered as  divisors on $Y$) can be naturally contracted to $\mathbf{CP}^1$ along mutually different directions. 
Let $\mu: Y\ra Z$ be this contraction and put $C_0=\mu(E)$, $\ol{C}_0=\mu(\ol{E})$.
Then the normal bundle of $C_0$ and $\ol{C}_0$ in $Z$ is $\mathscr O(1-n)^{\oplus 2}$.
This $Z$ is the twistor space of a LeBrun metric with torus action.

Finally a $\mathbf C^*\times\mathbf C^*$-action on the twistor space $Z$ has to be introduced.
On $\mathbf{P}(\mathscr E)$ it is explicitly given by
\begin{equation}\label{eqn-action}
(u,v,x,y,z)\mapsto (su,sv,tx,s^nt^{-1}y,z),\quad (s,t)\in\mathbf C^*\times\mathbf C^*,
\end{equation}
which preserves $X$ and fixes $p_0$ and $\ol{p}_0$.
When restricted to $U(1)\times U(1)$ this action commutes with the real structure (\ref{eqn-rs}).

\vsp\noindent
\textbf{(2.2)}
In the sequel we write $G=\mathbf C^*\times\mathbf C^*=\{(s,t)\}$ for simplicity.
To calculate $G$-action on $H^1(\Theta_Z)$, we introduce various $G$-equivariant exact sequences related to this cohomology group.
Our calculation in this subsection is similar to that of LeBrun in \cite{LB92} with some simplifications.
We note that the dimensions of the cohomology groups $H^i(\Theta_Z)$ are different from LeBrun's case in \cite{LB92} for $i=0,1$.

Let $\pi:Y\ra Q$ be the projection which is the composition of the small resolution $Y\ra X$ and the projection $X\ra Q$.
We have the following exact sequence of sheaves of $\mathscr O_Y$-modules
\begin{equation}\label{eqn-4terms}
0\lra \Theta_{Y/Q}\lra\Theta_Y\lra\pi^*\Theta_Q\lra \mathscr G\lra 0,
\end{equation}
where $\Theta_{Y/Q}$ and $\mathscr G$ denote the kernel and the cokernel of the natural homomorphism $\Theta_Y\lra \pi^*\Theta_Q$ respectively.
We decompose \eqref{eqn-4terms} into the following two short exact sequences:
\begin{equation}\label{eqn-exs1}
0\lra\Theta_{Y/Q}\lra\Theta_Y\lra \mathscr F\lra 0,
\end{equation}
\begin{equation}\label{eqn-exs2}
0\lra \mathscr F\lra \pi^*\Theta_Q\lra \mathscr G\lra 0,
\end{equation}
where $\mathscr F$ denotes the image sheaf of $\Theta_Y\lra \pi^*\Theta_Q$.
On the other hand we have a natural isomorphism $\Theta_{Y/ Q}\simeq \mathscr O_Y(E+\ol{E})$ and an exact sequence
\begin{equation}\label{eqn-exs3}
0\lra \mathscr O_Y\lra \mathscr O_Y(E+\ol{E})\lra\mathscr O_E(E)\oplus\mathscr O_{\ol{E}}(\ol{E})\lra 0.
\end{equation}
As are already explained we have $\mathscr O_E(E)\simeq\mathscr O_E(-1,1-n)$ and $\mathscr O_{\ol{E}}(\ol{E})\simeq\mathscr O_{\ol{E}}(1-n,-1)$.
By taking the direct image of  \eqref{eqn-exs3}, we obtain  an exact sequence
\begin{equation}\label{eqn-exs4}
0\lra \mathscr O_Q\lra\pi_*\mathscr O_Y(E+\ol E)\lra\mathscr O_Q(-1,1-n)\oplus\mathscr O_{Q}(1-n,-1)\lra 0,
\end{equation}
since $R^1\pi_*\mathscr O_Y=0$.
Because the relevant extension group $H^1(\mathscr O_Q(1,n-1)\oplus\mathscr O_{Q}(n-1,1))$ vanishes, \eqref{eqn-exs4} splits and we get $\pi_*\mathscr O_Y(E+\ol E)\simeq \mathscr O_Q\oplus\mathscr O_Q(-1,1-n)\oplus\mathscr O_{Q}(1-n,-1)$.
From this we obtain $H^i(\Theta_{Y/Q})\simeq H^i(\pi_*\mathscr O_Y(E+\ol{E}))\simeq H^i(\mathscr O_Q\oplus\mathscr O_Q(-1,1-n)\oplus\mathscr O_{Q}(1-n,-1))$, which vanishes if $i\ge 1$.
Therefore by \eqref{eqn-exs1} we obtain 
\begin{equation}\label{eqn-isom2}
H^i(\Theta_Y)\simeq H^i(\mathscr F)\quad\text{for }i\ge 1.
\end{equation}
On the other hand we have $H^i(Y, \mathscr G)\simeq \oplus _{i=1}^nH^i(C_i,N_{C_i/Q})\simeq \oplus_{i=1}^nH^i(\mathscr O_{C_i}(2))$ for any $i\ge0$.
Thus we obtain from (\ref{eqn-exs2}) an exact sequence
\begin{equation}\label{eqn-exs5}
0\lra H^0(\mathscr F)\lra H^0(\Theta_Q)\lra \oplus_{i=1}^nH^0(N_{C_i/Q})\lra H^1(\mathscr F)\lra  0
\end{equation}
and $H^i(\mathscr F)\simeq H^i(\pi^*\Theta_Q)\simeq H^i(\Theta_Q)=0$ for $i\ge 2$.
In particular, by (\ref{eqn-isom2}), we obtain 
\begin{equation}\label{eqn-isom3}
H^i(\Theta_Y)=0\quad \text{for }i\ge 2.
\end{equation}
Since any $C_i$ is a member of the pencil of $G$-invariant $(1,1)$-curves on $Q$, the image of the map $H^0(\Theta_Q)\lra\oplus_{i=1}^n H^0(N_{C_i/Q})$ in \eqref{eqn-exs5} is $6-1=5$-dimensional.
(This is more concretely shown in the proof of Proposition \ref{prop-main} below.)
It follows that $H^1(\Theta_Y)$ is $(3n-5)$-dimensional.

Associated to the blowing-down map $\mu:Y\ra Z$ we have a natural isomorphism 
$$
\Theta_{Y,E+\ol{E}}\simeq \mu^*\Theta_{Z,C_0+\ol{C}_0},
$$
where for a complex manifold $A$ and its complex submanifold $B$, $\Theta_{A,B}$ denotes the sheaf of holomorphic vector fields on $A$ which are tangent to $B$ in general.
On the other hand we readily have $H^i(\Theta_{Y,E+\ol{E}})\simeq H^i(\Theta_Y)$ for any $i\ge 0$ and $H^i(\mu^*\Theta_{Z,C_0+\ol{C}_0})\simeq H^i(\Theta_{Z,C_0+\ol{C}_0})$ for any $i\ge 0$.
Consequently we obtain a natural isomorphism 
\begin{equation}\label{eqn-isom4}
H^i(\Theta_Y)\simeq H^i(\Theta_{Z,C_0+\ol{C}_0}) \quad \text{for any }i\ge 0.
\end{equation}
Hence by \eqref{eqn-isom3} we obtain $H^i(\Theta_{Z,C_0+\ol{C}_0})=0$ for $i\ge 2$.
Therefore by an obvious exact sequence
\begin{equation}
0\lra\Theta_{Z,C_0+\ol{C}_0}\lra\Theta_Z\lra N_{C_0/Z}\oplus N_{\ol C_0/Z}\lra 0
\end{equation}
and \eqref{eqn-isom4} and $N_{C_0/Z}\simeq\mathscr O(1-n)^{\oplus 2}\simeq N_{\ol{C}_0/Z}$, we get
an exact sequence
\begin{equation}\label{eqn-exs6}
0\lra H^1(\Theta_Y)\lra H^1(\Theta_Z)\lra H^1(N_{C_0/Z})\oplus H^1(N_{\ol{C}_0/Z})\lra 0.
\end{equation}
It follows that the dimension of $H^1(\Theta_Z)$ is $(3n-5)+2\cdot 2(n-2)=7n-13$.
Also we obtain from the long exact sequence and \eqref{eqn-isom3}  that $H^2(\Theta_Z)=0$.

\vspace{3mm}\noindent
\textbf{(2.3)}
Now we have finished preliminaries for calculating the torus action on the cohomology group.
By the exact sequence (\ref{eqn-exs6}) which is obviously $G$-equivariant, it suffices
to calculate $G$-actions on $H^1(\Theta_Y)$ and $ H^1(N_{C_0/Z})\oplus H^1(N_{\ol{C}_0/Z})$ respectively.
To put the result in simple form, we use the following notation for expressing torus actions:
if a complex vector space $V$ of finite dimension $k$ is acted by the torus $G=\mathbf C^*\times \mathbf C^*=\{(s,t)\}$, $V$ can be decomposed essentially in a unique way into the direct sum of 1-dimensional $G$-invariant subspaces $V_i$, $1\le i\le k$.
For each $V_i$, $G$-action on $V_i$ takes the form $v_i\mapsto s^{m_i}t^{n_i}v_i$ for some integers $m_i$ and $n_i$.
Under this situation we write the $G$-action on $V$ by $\{(m_1,n_1),(m_2,n_2),\cdots,(m_k,n_k)\}$.
Then our result is as follows:

\begin{prop}\label{prop-main}
Let $Z$ be the twistor space of a LeBrun metric with torus action on $n\mathbf{CP}^2$, $n\ge 3$.
Then the natural action of the torus  on the cohomology group $H^1(Z,\Theta_Z)\simeq \mathbf C^{7n-13}$ is the direct sum of the following three representations of the torus:
\begin{equation}\label{eqn-rep1}
\big\{
\underbrace{(0,0),\cdots,(0,0)}_{n-1},
\underbrace{(1,0),\cdots,(1,0)}_{n-2},
\underbrace{(-1,0),\cdots,(-1,0)}_{n-2}
\big\}
\end{equation}
on $H^1(\Theta_Y)\simeq \mathbf C^{3n-5}$, and 
\begin{equation}\label{eqn-rep2}
\big\{
\underbrace{(1-n,1),(2-n,1),\cdots,(-2,1)}_{n-2},
\underbrace{(2-n,1),(3-n,1),\cdots,(-1,1)}_{n-2}
\big\}
\end{equation}
on $ H^1(N_{C_0/Z})\simeq\mathbf C^{2n-4}$, and 
\begin{equation}\label{eqn-rep3}
\big\{
\underbrace{(n-1,-1),(n-2,-1),\cdots,(2,-1)}_{n-2},
\underbrace{(n-2,-1),(n-3,-1),\cdots,(1,-1)}_{n-2}\big\}
\end{equation}
on $H^1(N_{\ol C_0/Z})\simeq\mathbf C^{2n-4}$.
\end{prop}


\noindent Proof.
First we prove that the torus action on $H^1(\Theta_Y)$ is as in \eqref{eqn-rep1}.
We use the exact sequence \eqref{eqn-exs5} which is also torus-equivariant sequence.
We first determine the image of the homomorphism $\alpha:H^0(\Theta_Q)\ra \oplus_{i=1}^n H^0(N_i)$ in \eqref{eqn-exs5}, where we write $N_i=N_{C_i/Q}$ for simplicity.
Viewing $H^0(\Theta_Q)$ as the Lie algebra of $\mathrm{Aut}_0(Q)\simeq \mathrm{PSL}(2,\mathbf C)\times \mathrm{PSL}(2,\mathbf C)$, $\alpha$ can be concretely given as follow:
for any $X\in sl(2,\mathbf C)\oplus sl(2,\mathbf C)$, let $\{A(t)\set t\in\mathbf C\}$ be the 1-parameter subgroup in $\mathrm{PSL}(2,\mathbf C)\times\mathrm{PSL}(2,\mathbf C)$ generated by $X$.
For any point $q\in C_i$, we associate the tangent vector at $q$ of the $A(t)$-orbit through $q$.
Consequently we obtain a tangent vector along $C_i$, which is a holomorphic section of $\Theta_Q|_{C_i}$.
Then projecting this onto $N_i$, we obtain an element of $H^0(N_i)$.
This is $\alpha(X)$.
In the sequel we choose a basis of $sl(2,\mathbf C)\oplus sl(2,\mathbf C)$ and for each member of the basis  we calculate their images under $\alpha$.

Before concretely calculating the image of $\alpha$, we give, for each $C_i$ ($1\le i\le n)$, a direct sum decomposition $\Theta_Q|_{C_i}\simeq \Theta_{C_i}\oplus N_i$ (namely, a splitting of 
$0\lra \Theta_{C_i}\lra\Theta_Q|_{C_i}\lra N_i\lra 0$).
For this, let $(u,v)$ be a non-homogeneous coordinate on $Q$ as in (2.1), and $\tau_i\in H^0(\Theta_{C_i})$ and $\nu_i\in H^0(\Theta_Q|_{C_i})$  holomorphic vector fields defined by
$$
\tau_i=\frac{\ptl}{\ptl u}+a_i\frac{\ptl}{\ptl v},\quad\nu_i=a_i\frac{\ptl}{\ptl u}-\frac{\ptl}{\ptl v}.
$$
Because $a_i$ is real, $\tau_i$ and  $\nu_i$ cannot be parallel and $\nu_i$ can be regarded as a (holomorphic) non-zero section of $\nu_i$.
Then we obtain a direct sum decomposition $\Theta_Q|_{C_i}\simeq \Theta_{C_i}\oplus N_i$.
Explicitly, if $\gamma=g(\ptl/\ptl u)+h(\ptl/\ptl v)$ is a holomorphic section of $\Theta_Q|_{C_i}$, we have 
\begin{equation}\label{eqn-dcp}
\gamma=\alpha\tau_i+\beta\nu_i;\quad \alpha=\frac{g+a_ih}{1+a_i^2},\quad
\beta=\frac{a_ig-h}{1+a_i^2}.
\end{equation}
Moreover, we can take $\{\nu_i,u\nu_i,u^2\nu_i\}$ as a basis of $H^0(N_i)$.

As a basis of $sl(2,\mathbf C)$ we choose
$$
A=
\begin{pmatrix}
1&0\\0&-1
\end{pmatrix},
\quad B=
\begin{pmatrix}
0&1\\0&0
\end{pmatrix},
\quad
C=
\begin{pmatrix}
0&0\\1&0
\end{pmatrix}.
$$
Corresponding 1-parameter subgroups are
$$
\begin{pmatrix}
e^t&0\\0&e^{-t}
\end{pmatrix},\quad
\begin{pmatrix}
1&t\\0&1
\end{pmatrix},\quad
\begin{pmatrix}
1&0\\t&1
\end{pmatrix}
$$
respectively, where $t\in\mathbf C$.
Then if we choose as a basis of $sl(2,\mathbf C)\oplus sl(2,\mathbf C)$ 
\begin{equation}\label{eqn-basis1}
(A,O),\; (B,O),\; (C,O),\; (O,A),\; (O,B),\; (O,C),
\end{equation}
where $O$ is the zero matrix, and if $\gamma_{i1},\gamma_{i2},\cdots,\gamma_{i6}\in H^0(N_i)$ denotes the image of the above 6 generators of  $sl(2,\mathbf C)\oplus sl(2,\mathbf C)\simeq H^0(\Theta_Q)$ by the homomorphism $H^0(\Theta_Q)\ra H^0(N_i)$ respectively, then we obtain by using \eqref{eqn-dcp}
\begin{multline}\label{eqn-image2}
\gamma_{i1}=\frac{a_i}{1+a_i^2}u\nu_i,\quad
\gamma_{i2}=-\frac{a_i}{1+a_i^2}u^2\nu_i,\quad
\gamma_{i3}=\frac{a_i}{1+a_i^2}\nu_i,\\
\gamma_{i4}=-\frac{a_i}{1+a_i^2}u\nu_i,\quad
\gamma_{i5}=\frac{a_i^2}{1+a_i^2}u^2\nu_i,\quad
\gamma_{i6}=-\frac{1}{1+a_i^2}\nu_i.
\end{multline}
Thus the image of each member of \eqref{eqn-basis1} by $\alpha$ is
$
\gamma_k:=\sum_{i=1}^n\gamma_{ik}\in\oplus_{i=1}^nH^0(N_i), \, 1\le k\leq 6
$
respectively.
Obviously $\gamma_1=-\gamma_4$ and it is easily verified that $\gamma_2,\gamma_3,\gamma_4,\gamma_5,\gamma_6$ are linearly independent (in $\oplus H^0(N_i)$).
Thus we have obtained 
\begin{equation}\label{eqn-basis5}
\mathrm{Image}(\alpha)=\langle\gamma_2,\gamma_3,\gamma_4,\gamma_5,\gamma_6\rangle\subset\oplus_{i=1}^nH^0(N_i).
\end{equation}

Now  we are able to calculate $G$-action on $H^0(N_i)$. 
Recall that by \eqref{eqn-action} we have $(u,v)\mapsto (su,sv)$ for $(u,v)\in Q$.
(In particular a subgroup $\{(s,t)\in G\set s=1\}$ acts trivially on $Q$.)
It follows that 
\begin{equation}\label{eqn-act7}
\nu_i\mapsto s\nu_i,\; u\nu_i\mapsto u\nu_i,\; u^2\nu_i\mapsto s^{-1}u^2\nu_i \quad\text{for}\quad s\in\mathbf C^*
\end{equation}
 for each of the basis of $H^0(N_i)$.
The $G$-action on $\oplus H^0(N_i)$ is the direct sum of these $n$ representations.
Needless to say, $\{\nu_i,u\nu_i,u^2\nu_i\set 1\le i\le n\}$ is a basis of  $\oplus_{i=1}^n H^0(N_i)$.
Instead of this basis,  it is easily seen by carefully looking \eqref{eqn-image2} that  we can take, as a basis of $\oplus_{i=1}^n H^0(N_i)$,
$$
\left\{\gamma_i,\,\nu_j,\,u\nu_k,\,u^2\nu_l\set 2\leq i\le 6,\,3\le j\le n,\,2\le k\le n,\,3\le l\le n\right\}.
$$
Combining this with \eqref{eqn-basis5}, we obtain 
\begin{multline}
\left(\oplus_{i=1}^nH^0(N_i)\right)/\,\mathrm{Image}\{\alpha:H^0(\Theta_{Q})\ra \oplus_{i=1}^nH^0(N_i)\}\\
\simeq
\left\{\nu_j,\,u\nu_k,\,u^2\nu_l\set 3\le j\le n,\,2\le k\le n,\,3\le l\le n\right\}.
\end{multline}
Hence by \eqref{eqn-exs6} we have obtained that the $G$-action on $H^1(\mathscr F)\simeq H^1(\Theta_Y)$ is given by (\ref{eqn-rep1}).

Our next task is to calculate $G$-action on $H^0(N_{C_0/Z}$).
For this, we first consider the following two  divisors
$$
D_0:=\{u=0\}\cap X\subset\mathbf P(\mathscr E),
\quad D_{\infty}=\{u=\infty\}\cap X\subset\mathbf P(\mathscr E)
$$
in $X$, which are clearly $G$-invariant.
Obviously $D_0$ and $D_{\infty}$ are disjoint.
If we use the same symbols to denote the corresponding $G$-invariant divisors in $Y$ and $Z$, $D_0\subset Z$ and $D_{\infty}\subset Z$ intersect transversally along $C_0$.
(Note that by the blowing-down $\mu:Y\ra Z$ the divisor $E$ is blown-down along fibers of the projection to the second factor of $E\simeq Q\simeq\mathbf{CP}^1\times\mathbf{CP}^1$.
On the other hand we do not need to be careful for the small resolution $Y\ra X$ since $E$ and $\ol{E}$ are disjoint from the singular points of $X$.)
Therefore by setting $\Gamma_0=D_0\cap E\subset X$ and $\Gamma_{\infty}=D_{\infty}\cap E\subset X$, we have
\begin{equation}
N_{C_0/Z}\simeq N_{C_0/D_0}\oplus N_{C_0/D_{\infty}}\simeq N_{\Gamma_0/D_0}\oplus N_{\Gamma_{\infty}/D_{\infty}}.
\end{equation}
Moreover we have $N_{\Gamma_0/D_0}\simeq \mathscr O(1-n)\simeq N_{\Gamma_{\infty}/D_{\infty}}$ since $N_{E/X}\simeq \mathscr O(-1,1-n)$.
Thus it suffices to determine $G$-actions on $H^1(N_{\Gamma_0/D_0})$ and $H^1(N_{\Gamma_{\infty}/D_{\infty}})$ respectively.
For these, we use \v{C}ech  representation of elements of $H^1(\mathscr O(1-n))$.
First we calculate $G$-action on $H^1(N_{\Gamma_0/D_0})$.
The point $(u,v,x,y,z)=(0,0,0,1,0)\in \mathbf P(\mathscr E)$ lies on $\Gamma_0$ and is a $G$-fixed point.
We use $v$ as a non-homogeneous coordinate on $\Gamma_0$.
Then over $\mathbf C\subset\Gamma_0$ on which $v$ is valid, one can use $z/y$ as a fiber coordinate of $N_{\Gamma_0/D_0}$.
Then by \eqref{eqn-action} $G$-action on the total space of $N_{\Gamma_0/D_0}$ is given by
$$
(v,(z/y))\mapsto (sv,s^{-n}t\,(z/y)),\quad (s,t)\in G.
$$
On the other hand, any element of $H^1(\mathscr O(1-n))$ is represented by a linear combination of the following $n-2$ sections of $\mathscr O(1-n)$ over $\mathbf C^*$;
$$
\zeta_k:v\mapsto v^{-k},\; v\in \mathbf C^*;\;1\leq k\leq n-2.
$$
Then since $s^{-n}t\cdot v^{-k}=s^{k-n}t\cdot(sv)^{-k}$, $\zeta_k$ is mapped to $s^{k-n}t\cdot\zeta_k$ by $(s,t)\in G$.
Thus in the notation we have introduced before Proposition \ref{prop-main}, we obtain that the $G$-action on $H^1(N_{C_0/D_0})\simeq H^1(N_{\Gamma_0/D_0})$ is given by 
\begin{equation}\label{eqn-act6}
\{(k-n,1)\set k=1,2,\cdots, n-2\}.
\end{equation}
Next we calculate $G$-action on $H^1(N_{\Gamma_{\infty}/D_{\infty}})$ in a similar way.
As a $G$-fixed point on $\Gamma_{\infty}$ we choose a point $(u,v,x,y,z)=(\infty,0,0,1,0)$ and again use $v$ as a non-homogeneous coordinate on $\Gamma_{\infty}$.
Then as a fiber coordinate of $N_{\Gamma_{\infty}/D_{\infty}}$ we can use $z/(u^{-1}y)$.
(The multiplication of $u^{-1}$ comes from $y\in \mathscr O(1,n-1)$.)
Again by \eqref{eqn-action},  in this coordinate the $G$-action on the total space of $N_{\Gamma_{\infty}/D_{\infty}}$ is given by
$$
\left(v, z/(u^{-1}y)\right)\mapsto \left(sv, \,s^{1-n}t\cdot(z/(u^{-1}y))\right).
$$
Since $s^{1-n}t\cdot u^{-k}=s^{k-n+1}t\cdot(su)^{-k}$ this time, we have that $\zeta_k$ is multiplied by $s^{k-n+1}t$ by $(s,t)\in G$.
It follows that $G$-action on $H^1(N_{C_{\infty}/D_{\infty}})\simeq H^1(N_{\Gamma_\infty/D_{\infty}})$ is given by 
\begin{equation}\label{eqn-act8}
\{(k-n+1,1)\set k=1,2,\cdots,n-2\}.
\end{equation}
By \eqref{eqn-act6} and \eqref{eqn-act8}, we obtain that the $G$-action on $H^1(N_{C_0/Z})$ is as in \eqref{eqn-rep2}.

Finally, the $G$-action on $H^1(N_{C_{\infty}/Z})$ is known to be given by \eqref{eqn-rep3} by taking $\ol{D}_0=\{v=0\}\cap X$ and $\ol{D}_{\infty}=\{v=\infty\}\cap X$ instead of $D_0$ and $D_{\infty}$ in the above argument.
\proofend

\vsp
The statement of Proposition \ref{prop-main} and its proof perfectly work also for the case $n=1$ and $n=2$ but in these cases it brings not much informations.

\section{Equivariant deformations of the metric and examples}
\noindent
\textbf{(3.1)}
Proposition \ref{prop-main} is not so useful in itself.
In this subsection, by using Proposition \ref{prop-main},  we give a geometric characterization of $U(1)$-subgroups for which there exists a $U(1)$-equivariant deformation which does not preserve full torus symmetry.
Let $E_1+E_2+\cdots+ E_{n-1}\subset Y$ be the exceptional curve of the small resolution of $p_0\in X$ given in (2.1), where
 $E_i\simeq \mathbf{CP}^1$ is the exceptional curve obtained in the $i$-th blow-up (along the 2-dimensional center we have explicitly given), so that $E_i$ and $E_{j}$ $(i\neq j)$ intersect and iff $|i-j|=1$.
Because any $E_i$ is not affected by the blowing-down $\mu:Y\ra Z$ we use the same notation to represent the corresponding rational curves in  $Z$.
Clearly $C_0$ and $\ol{C}_0$ are disjoint from $E_1+\cdots+ E_{n-1}\subset Z$.
The curve $\{y=u=v=0\}$ in $X$ connects $p_0$ and $\ol{E}$.
Let $B_0\subset Z$ be the strict transform of this curve.
$B_0$ connects $\ol{C}_0$ and $E_1$.
Similarly the rational curve $\{x=u=v=0\}\subset X$ connects $p_0$ and $E$, and its strict transform in $Z$ is denoted by $B_n$ which connects $E_{n-1}$ and $C_0$.
In this way we obtain a string of $(n+3)$ smooth rational curves
\begin{equation}\label{eqn-string}
\ol{C}_0+B_0+E_1+E_2+\cdots+E_{n-1}+B_n+ C_0,
\end{equation}
where only adjacent two curves intersect.
Adding the conjugate curves $\ol{B}_0+\ol{E}_1+\cdots+\ol{E}_{n-1}+\ol{B}_n$ to \eqref{eqn-string}, we obtain a cycle of $(2n+4)$ rational curves in $Z$.
Obviously this cycle of rational curves are $G$-invariant and the intersection points of the irreducible components are (isolated) $G$-fixed points of $Z$.
Moreover, this cycle is the basel locus of the pencil of $G$-invariant divisors in $|(-1/2)K_Z|$.
Note that the image of this cycle onto $n\mathbf{CP}^2$ by the twistor fibration is a cycle of torus invariant $(n+2)$ spheres, on which some of $U(1)$-subgroup of the torus acts trivially.

Elements of the torus $U(1)\times U(1)\subset G$  fixing any point of $C_0$ form a $U(1)$-subgroup, which we denote by $K_0$.
By reality, $K_0$ automatically fixes any point of $\ol C_0$.
Similarly let $K_i\subset U(1)\times U(1)$, $1\le i\le n-1$, be the $U(1)$-subgroup fixing any point of $E_i$ (and hence $\ol E_i$).
In this way we have obtained $n$ numbers of $U(1)$-subgroups in the torus
(so that in particular we do not consider $U(1)$-subgroup fixing $B_0$ and $B_n$ among the cycle above).

\begin{prop}\label{prop-moduli}
Let $K$ be any $U(1)$-subgroup in the torus.
Then  LeBrun's metric with torus action on $n\mathbf{CP}^2$, $n\ge 3$, can be $K$-equivariantly deformed into self-dual metric with only $K$-isometry  if and only if $K=K_i$ for some $i$, $0\le i\le n-1$.
Moreover, the dimension of the moduli spaces of resulting self-dual metrics with just $U(1)$-isometry obtained in this way become as follows:
\begin{itemize}
\item
$(3n-6)$\text{-dimensional for $K_0$-equivariant deformations},

\item
$n$\text{-dimensional for $K_i$-equivariant deformations for $i=1$ or $n-1$},

\item
$(n+2)$\text{-dimensional for $K_i$-equivariant deformations for $2\leq i\le n-2$}.
\end{itemize}
Furthermore, in the second and the third cases, the self-dual metric is not conformally isometric to LeBrun metric.
(Note that if $n=3$ the third item does not occur.)
\end{prop}

\noindent
Proof. 
The $G$-action on $C_0$ and the exceptional curves $E_i$ $(1\le i\le n-1)$ can be readily computed by using \eqref{eqn-action} and explicit small resolution given in (2.1).
Consequently we obtain that the subgroups $K_i$ are explicitly given by
$$
K_0=\{(s,t)\in U(1)\times U(1)\set s=1\},
$$
$$
K_i=\{(s,t)\in U(1)\times U(1)\set t=s^i\},\;1\le i\le n-1.
$$
Then comparing these with the result in Proposition \ref{prop-main}, we obtain that for a $U(1)$-subgroup $K\subset U(1)\times U(1)$, the $K$-fixed subspace $H^1(\Theta_Z)^{K}$ contains $H^1(\Theta_Z)^{U(1)\times U(1)}$ as a proper subspace if and only if $K=K_i$ for some $i$, $0\le i\le n-1$.
Noting that $H^1(\Theta_Z)^{K}$ is the tangent space of the Kuranishi family of  $K$-equivariant deformations of $Z$ (since $H^2(\Theta_Z)=0$),
it follows that $Z$ admits a $K$-equivariant deformation which does not preserve the full torus symmetry if and only if $K=K_i$ for some $i$, $0\le i\le n-1$.
Since the $U(1)\times U(1)$-action on $H^1(\Theta_Z)$ commutes with the natural real structure induced by that on $Z$, the situation remains unchanged even after restricting to the real part of $H^1(\Theta_Z)$;  
namely $Z$ admits a $K$-equivariant deformation which preserves the real structure but does not preserve the full torus symmetry if and only if $K=K_i$ for some $i$, $0\le i\le n-1$.
This implies that LeBrun's twistor space admits a non-torus equivariant, $K$-equivariant deformation as a twistor space if and only if $K=K_i$ for some $0\le i\le n-1$.
Going down on $n\mathbf{CP}^2$, we obtain the first claim of the proposition.

Next we compute the dimension of the moduli space by using Proposition \ref{prop-main}.
For $K_0$-equivariant deformation, we obtain from \eqref{eqn-rep1}--\eqref{eqn-rep3} that $H^1(\Theta_Z)^{K_0}$ is just $H^1(\mathscr F)$ that is  $(3n-5)$-dimensional.
On this subspace the quotient torus $(U(1)\times U(1))/K_0$ acts non-trivially and its orbit space is just the (local) moduli space of $K_0$-equivariant self-dual metrics on $n\mathbf{CP}^2$. 
In particular its dimension is $(3n-5)-1=3n-6$.
For $K_1$ and $K_{n-1}$-equivariant deformations, the fixed subspace $H^1(\Theta_Z)^{K_i}$ is ($(n-1)+2=n+1$)-dimensional.
Therefore the moduli space is $n$-dimensional.
For other equivariant deformations, we have $H^1(\Theta_Z)^{K_i}$, $2\le i\le n-2$, is ($(n-1)+2\cdot2=n+3)$-dimensional and the  moduli space becomes $(n+2)$-dimensional.

Finally it is easily seen that the action of $K_i=\{(s,t)\set t=s^i\}$, $1\le i\le n-1$,
 on the torus-invariant rational curve $B_0$ is explicitly given by $\tilde x\mapsto s^i\tilde x$ for an affine coordinate $\tilde x$ on $B_0$.
This means that if $i\ge 2$ then $K_i$ contains non-trivial isotropy along $B_0$.
Therefore by a theorem of LeBrun \cite{LB93} characterizing LeBrun metric by semi-freeness of the $U(1)$-action, we conclude that self-dual metric obtained by $K_i$-equivariant, non-torus equivariant deformation is not conformally isometric to LeBrun metric.
For the remaining $K_1$-equivariant deformation, it suffices to consider $B_n$ instead of $B_0$.
\proofend


\vsp
\noindent \textbf{(3.2)}
Finally we discuss some examples.
\begin{example}{\rm{
First we consider torus equivariant deformation of LeBrun's metric with torus action on $n\mathbf{CP}^2$.
By Proposition \ref{prop-main} the subspace of $H^1(\Theta_Z)$  consisting of  vectors which are torus-invariant is $(n-1)$-dimensional.
This is consistent with the fact that the moduli space of  LeBrun's metrics with torus action (or more generally, Joyce's metric with torus action \cite{J95}) is $(n-1)$-dimensional.
See also a work of Pedersen-Poon \cite{PP95}, where the dimension of the moduli space is calculated via a construction of Donaldson and Friedman.
}}
\end{example}

\begin{example}{\rm{
Consider $K_0$-equivariant deformation of LeBrun's metric with torus action.
By definition $K_0$ fixes any point of $C_0$ and $\ol{C}_0$ and acts semi-freely on $n\mathbf{CP}^2$.
By Proposition \ref{prop-moduli} the moduli space of self-dual metrics on $n\mathbf{CP}^2$ obtained by $K_0$-equivariant deformation is $(3n-6)$-dimensional.
Of course this coincides with the moduli number obtained by LeBrun \cite{LB91, LB92}.
(LeBrun's result is much stronger in that  his construction makes it possible to determine the  {\em global structure} of the  moduli space.)}}
\end{example}

\begin{example}{\rm{
Let $n=3$ and 
consider $K_1$-equivariant deformation of LeBrun's metric with torus action on $3\mathbf{CP}^2$.
By Proposition \ref{prop-moduli} the moduli space of self-dual metrics on $3\mathbf{CP}^2$ obtained by $K_1$-equivariant deformation of LeBrun metrics with torus action is $3$-dimensional.
Since $K_1$ does not act semi-freely on $3\mathbf{CP}^2$, these self-dual metrics are not conformally isometric to the LeBrun metric (obtained by so called hyperbolic ansatz).
In a recent paper \cite{H-bitan} the author determined a global structure of this moduli space.
In particular, the moduli space is connected and $3$-dimensional, which is equal to the dimension obtained in Proposition \ref{prop-moduli}.
We note that the situation for $K_2$-equivariant deformations is completely the same, since $K_1$-action and $K_2$-action are interchanged by a diffeomorphism of $3\mathbf{CP}^2$.
This is always true for $K_1$-action and $K_{n-1}$-action for any $n\,(\ge 3)$.
It is also possible to show that the twistor space obtained by $K_1$-equivariant deformations of LeBrun metric with torus action on $n\mathbf{CP}^2$  is, at least for small deformations, always {\em Moishezon}.
 }}
\end{example}

\begin{example}{\rm{
In \cite{Hon01} it was prove that being {\em Moishezon} twistor space is not preserved under $\mathbf C^*$-equivariant small deformations as a twistor space.
This is obtained by letting $n=4$ and considering $K_2$-equivariant small deformations of LeBrun twistor spaces with torus action. 
This in particular implies that if one drops the assumption of the semi-freeness of $U(1)$-isometry, then the twistor space is not Moishezon in general.
}}
\end{example}

\small
\vspace{13mm}
\hspace{5.5cm}
$\begin{array}{l}
\mbox{Department of Mathematics}\\
\mbox{Graduate School of Science and Engineering}\\
\mbox{Tokyo Institute of Technology}\\
\mbox{2-12-1, O-okayama, Meguro, 152-8551, JAPAN}\\
\mbox{{\tt {honda@math.titech.ac.jp}}}
\end{array}$

\end{document}